\documentclass{siamltex}		
\usepackage[utf8]{inputenc}
\usepackage{graphicx,multirow}  
\usepackage{verbatim,float} 
\usepackage{todonotes}
\usepackage{amssymb,amsmath}
\usepackage{graphicx,amsmath,amsfonts,amssymb,subfigure}
\usepackage{color}
\usepackage[ruled,vlined,linesnumbered]{algorithm2e}

\def\RR{\mathcal{R}}

\def\nref#1{(\ref{#1})}
\def\inv{^{-1}}%

\def\Sup#1{^{(#1)}}%
 
\newcommand{\eq}[1]{\begin{equation}\label{#1}}
\newcommand{\en}{\end{equation}}
\graphicspath{{../talk/FIGS/},{../../DMR_EPS/},{./../../Tutor/EPS/}}

\title{
  Iterative methods for linear systems of equations:
 A brief historical journey}

\author{Yousef Saad
\thanks{Address: Computer Science \& Engineering,
  University of Minnesota, Twin Cities.
  Work supported by NSF grant DMS 1521573.
  e-mail:
  \texttt{ saad@umn.edu}
  }
}
\begin{document} 

\maketitle 

\begin{abstract}
  This  paper   presents  a brief historical survey of iterative methods
  for solving linear systems of equations. The journey begins with Gauss
  who developed the first 
  known method that can be termed iterative.
  The early 20th century saw good progress of these methods
  which were initially used to solve least-squares systems, and then
  linear systems arising from the discretization of partial different
  equations. Then iterative methods received
  a big impetus in the  1950s - partly because of the development of
  computers. The survey does not attempt to be exaustive. Rather, the aim is
  to underline the way of thinking at a specific time and to highlight the
  major ideas that steered the field.
\end{abstract}


\section{It all started with Gauss}\label{sec:Gauss} 
A big part  of the contributions of Carl Friedrich  Gauss can be found
in the voluminous exchanges he  had with contemporary scientists. This
correspondance has  been preserved in a  number of books, e.g., in the
twelve `Werke' volumes  gathered from 1863 to 1929 at
the                            University                           of
G\"ottingen\footnote{https://gdz.sub.uni-goettingen.de/volumes/id/PPN235957348}.
There are also books specialized  on specific correspondances.  One of
these  is dedicated  to the  exchanges  he had  with Christian  Ludwig
Gerling \cite{GaussGerling}.   Gerling was a student of Gauss  under
whose supervision he  earned a PhD from the  university of G\"ottingen
in 1812. Gerling  later became professor of  mathematics, physics, and
astronomy at the  University of Marburg where he spent the rest of his
life from  1817 and
maintained a  relation  with Gauss  until
Gauss's death in  1855.  We learn from  \cite{SmithRev1928} that there
there were  388 letters exchanged  between the  two from 1810  to 1854
(163 written by Gauss and 225 by Gerling).

It  is  one of  these  letters  that  Gauss  discusses his  method  of
\emph{indirect  elimination} which  he contrasted  with the  method of
\emph{direct   elimination}  or   Gaussian   elimination  in   today's
terminology.  Gauss wrote  this letter to Gerling on  Dec. 26th, 1823.
Gerling was a specialist of geodetics  and the exchange in this letter
was  about the  application of  the method  of least-squares  which
Gauss invented in the early 1800s to geodetics.
An English translation of this letter was published in  1951 by George Forsythe
\cite{Letter2Gerling}.
The historical context of this translation is interesting in itself.
In the forward, Forsythe begins  by stating that his specific aim was
to find whether or not the
reference in a book by Zurm\"uhl on matrices
\cite{Zurmul50} that mentioned  the method of relaxation
described by Southwell \cite{Southwell-1935,SouthwellBook} is the same as the one given in
``Dedekind's report on Gauss's Lectures''
(see references in \cite{Letter2Gerling})
\emph{``...It is believed by some computers\footnote{
  In the 1950s a `computer' often referred to someone who specialised
  in  numerical computing, i.e., 
  whose job was to carry out a calculation given a certain algorithm,
  or set of formulas.}
  however, that Gauss's method is a
different one, namely the related method given by Seidel
\cite{Seidel1874}.
In the interest of giving Gauss his proper credit as a proponent of
relaxation, the following translation of a letter by Gauss is offered''.
''}

Let us take a look at the content of the letter.
After some manipulations of data, Gauss arrives 
at a least-squares system for determining angles formed
by certain directions for 4 German cities. For this he states that he will
describe an `indirect' method for solving the normal equations.
The least-squares system is stated as follows:
\begin{alignat}{5}
  0 &  =   + 6   & &  + 67a  & & -13b   & & -28c   & &- 26d \\
  0 &  =  - 7558 & &-13a   & &+ 69b  & &- 50c   & &- 6d\\
  0 &  =  -14604 & &-28a   & &- 50b  & &+ 156c  & &- 78d\\
  0 &  =+ 22156  & &-26a   & &- 6b   & &- 78c   & &+ 110 d
\end{alignat}
As can be observed the sum of the right-hand side is zero (the right-hand side
entries are termed `absolute terms' by Gauss). In fact all columns add up to zero.
The preliminary step to produce a system of this type is a
`trick' which helps get a better convergence and provides what we may call today
invariant, i.e., a property that is always satisfied and can therefore be useful in detecting
calculation errors.

In the iterative scheme  proposed, the coordinates of the
solution change and the right-hand side, which is in fact the residual
vector, is updated each time. The mechanical nature of the procedure is
at the same time simple and appealing.
So for example, all coordinates of the solution are set to zero, and in the first step, he selects to modify the
4th coordinate because it would lead to the biggest decrease in the residual.
The letter then shows the progress of the algorithm in the following table:
\[
\begin{array}{r|r|r|r|r|r|r|r|}
      &  d=-201  & b=+92  & a=-60  & c=+12  & a=+5  & b=-2 & a=-1\\ \hline
+ 6   & +5232    & +4036  & + 16  & -320    & + 15  & +41  & -26\\
- 7558 &  -6352  & - 4    & + 776  & +176   & +111  & -27  & -14\\
-14604 &  +1074  & -3526  &  -1846 &  + 26  & -114  & -14  & + 14\\
+22156 & + 46    & - 506  & +1054  & +118   & - 12  & 0    & +26
\end{array}
\]
The first column is just the initial residual (corresponding to the initial guess $x=0$ in today's
terminology). The largest entry is the 4th, and so Gauss selects to update $d$ which now has the value $d = -201$
(obtained by making the last equation satisfied). Note that all values are rounded ans this is one of the important
attractions of this method. The second column shows the next modification to be added this time to $b$.
In today's notation we would write something like $\Delta b = + 92$.
At the end of the 7 steps above,
we have $a = -56, b = +90, c= 12, d = -201$.
The new residual is shown in the second column. The process is continued in this way.
It turns out for this example,  nothing changes after the 7 steps shown above: 
\emph{    Insofar as I carry the calculation only to the nearest 2000-th of a second,
  I see that now there is nothing more to correct...} and the final solution is displayed.
Gauss ends the letter with the words:
\emph{
...
  Almost every evening I make a new edition of the tableau, wherever
     there is easy improvement. Against the monotony of the surveying business,
     this is always a pleasant entertainment; one can also see immediately
     whether anything doubtful has crept in, what still remains to be desired,
etc. I recommend this method to you for imitation. You will hardly ever
again eliminate directly, at least not when you have more than 2.
The indirect procedure can be done while half asleep, or while thinking
about other things.
}

Gauss  recommends this iterative scheme \emph{(indirect elimination)}
over  Gaussian elimination for systems of order $>2$.
We will contrast this with other recommendations later.

This appears to be the first known  reference to a use of an iterative
method   for  solving   linear   systems.   Later,   in  1845   Jacobi
\cite{jacobi45} developed a relaxation type method in which the latest
modification is not immediatly incorporated  into the system.  In that
same paper he  introduces a way of modifying a  linear system by using
what we  now know as  `Jacobi rotations' to annihilate  large diagonal
entries before  performing relaxations. This is in order  to speed-up
convergence of the iteration, and so one can say that
Jacobi introduced the first known form of preconditioning.
That  same technique that uses cleverly
selected  rotations, was  advocated in  a paper  that appeared  a year
later \cite{jacobi46} for solving symmetric eigenvalue problems.

Though this letter to Gerling  dealt with a $4 \times 4$ linear system,
Gauss solved bigger  systems with the same method. 
Bodewig \cite[p. 145]{Bodewig-Book} describes the state of the art with this method
in its early stages as follows:
\emph{
     Gauss had systems of 20—30—40 unknowns, later still
higher systems were met, for instance, in the triangulation of Saxony
by NAGEL in 1890: 159 unknowns, H. BOLTZ at Potsdam 670 unknowns
and in the present triangulation of Central Europe several thousands
of unknowns.}  


\section{Solution by relaxation}
The term `relaxation' which is in common use today to describe the general process first invented by Gauss,
seems to have been introduced by Southwell \cite{Southwell-1935,SouthwellBook}.
Suppose we have a linear system of the 
form
\eq{eq:Ax=b}
A x = b .
\en
which can be viewed  as a collection of equations:
\eq{eq:Aix=bi}
a_i x = \beta_i, \quad i=1,\cdots, n
\en
where  $a_i  \in \ \RR^{1 \times n}$ is the $i$-th row of $A$ and $\beta_i$ the $i$-th component of $b$.
We will denote by $r$ the residual vector:
\eq{eq:res}
r = b -A x . \en
The relaxation process is then as follows where we set $r\Sup{new} =  b - A x\Sup{new}$: 
\begin{center}  
    Modify $i$-th component of $x$ into $x_i\Sup{new} :=x_i+\delta $ so that:
    $r_i\Sup{new}  = 0$.
  \end{center}
This means that we should have   
    \[
      a_i (x + \delta e_i) = \beta_i \quad \longrightarrow \quad \delta = \frac{r_i}{a_{ii} }
    \]
    This is done in a certain order. In the original approach by Gauss, the method is to select $i$ to be
    the coordinate of the residual that has the largest entry. The newly computed component of $x$ is replaced into
    $x$ and the new associated residual is then computed. In the Jacobi method, all $n$ components are updated using the
    same $r$.  
    
    In 1874, Seidel~\cite{Seidel1874} described a relaxation method that was again geared toward solving normal equations.
    His method can be viewed as a relaxation for the system $A^T A x = A^T b$ of normal equations,
    and because of this specificity he was able to argue     for convergence.
    He also mentions that unknowns need not be processed from 1 to $n$ cyclically. Instead
     convergence is improved  by making each dominant residual into
zero, which is the same scheme as the one initially proposed by Gauss.
    In the same paper Seidel also developed a block method whereby a few unknowns are  proccessed at the same time.
    Often in the literature that followed, `relaxation' became synonymous with Seidel's method and the method was often called
    Seidel's method. It is now called the Gauss-Seidel method in an effort to give credit to Gauss who invented the non-cyclic
    variant almost 50 years before him. The fact that  Seidel recommends against processing the unknown cyclically
    prompted Gorge Forsythe to remark that    
\emph{
``the Gauss-Seidel method was not known to Gauss and not                       
recommended by Seidel''} according to Householder,
see \cite[p. 115]{Householder}. In the same notes, 
Householder also mentions that  Nekrasov, a Russian author 
(see \cite{Householder} for the reference) defined the exact same method 
as Seidel in 1884 and that ``the method is called Nekrasov's method in the
Russian literature''. In fact, to this day the method is referred to as
the method of Seidel-Nekrasov by some Russian authors.
Nekrasov analyzed the method theoretically
\cite{Nekrasov1885} and the paper
\cite{MehmkeNekrasov92} shows a convergence result.

One of the main attractions of the cyclic version of Gauss-Seidel iteration is that it can easily
programmed or `mechanized' as was said in the early days of computing.
David Young recounts the following anecdote in \cite{YoungHist}: 
\emph{
``Not too long after I began my work, Sir
Richard Southwell visited Birkhoff at Harvard.
One day when he, Birkhoff and I were together, I
told him what I was trying to do. As near as I
can recall, his words were
\emph{``any attempt to mechanize relaxation
methods would be a waste of time."}
This was somewhat discouraging, but my propensity
of making numerical errors was so strong that I
knew that I would never be able to solve significant problems except by machines. Thus, though 
discouraged, I continued to work.''}

Relaxation-type methods can be written in the form of fixed point
iterations and this makes it easy to analyze their convergence.
Consider the  decomposition $ A = D - E - F $ where $D$ is a diagonal matrix
having the same (diagonal) entries as those of $A$,
$-E$ is the stricly lower triangular part of $A$ and
$-F$ its stricly upper triangular part.
Then the method of Gauss-Seidel generates the iterates defined by
\eq{eq:GS}  x^{(k+1)} = (D-E) \inv \left(F x^{(k)}  + b \right) .
\en
Here one can write $Ax = b$ as $(D-E) x = F +b $, using the
\emph{splitting} $A = (D-E) \ - \ F$, from which the above iteration follows.
Similarly, the  Jacobi iteration is of the form
$  x^{(k+1)}  = D\inv [(E+F) x^{(k)}  + b ] $
and is based on the splitting: $A =  D \ - \ (E+F) $.

In addition to the original paper by Seidel mentioned above,
convergence of the standard Gauss-Seidel process was studied early on by several authors,
see, e.g., \cite{vonMisesGeiringer}. A number of these results can found in a chapter of
Bodewig's book \cite[Chap. 7]{Bodewig-Book} that 
starts with the warning:
\emph{But, first, let us note that these theorems are more or less superfluous
in practical computation. For the iteration methods will only be applied
when the convergence is evident at first sight, that is, when the diagonal
dominates strongly whereas in other cases convergence will be too slow
even for modern computing machines so that it is better to apply a direct
procedure (Gauss or variants).}

\section{Early 20th century}
The early 20th century was marked  by the beginning of the application
of  iterative  methods to  problems  modeled  by partial  differential
equations. Up to that period, iterative methods were mainly utilized to
solve linear systems that originated from normal equations.
On the other hand, a method proposed  by Liebmann~\cite{Liebmann}
was geared specifically toward solving discretized
Poisson equations.
The method is nothing but what we term today the Gauss-Seidel method,
and for this reason the Gauss-Seidel iteration when applied to Partial Differential Equations
was often called  the Liebmann method.  
It is known  as Nekrasov's method in the Russian literature
\cite{Nekrasov}.

In a  remarkable paper published in 1910,  Richardson \cite{Richardson-1910}
put together a number of  techniques for solving simple PDEs (Laplace,
Poisson, Bi-Harmonic, ..) by finite differences.  He then describes an
iterative  scheme for  solving  the linear  system  that results  from
discretizing these equations.  The PDEs addressed in the paper are all
of a homogenous type, e.g., Laplace,  or $\Delta u = 0$, with boundary
conditions. This  results in a  linear system  that can be  written as
$Ax + b =  0$ where $A$ acts on interior points  only and $b$ reflects
the action of  the discretized operator on the  boundary points. Thus,
for an arbitrary $x$, the vector $A  x + b$ represents the residual of
the system under  consideration.  With this notational  point in mind,
the method  introduced by Richardson in  this paper can be  written in
the form:  \eq{eq:Rich} x_{j+1} =  x_j - \frac{1}{\alpha_j} A  r_j \en
and results  in a polynomial  iteration scheme whose residual  at step
$k$ satisfies
\eq{eq:resk} r_k = (I  - \frac{A}{\alpha_k } )\cdots (I -
  \frac{A}{\alpha_2}) \cdots (I - \frac{A}{\alpha_1}) r_0 .  \en
  Thus, $r_k$  is of the form:
  \eq{eq:resk1} r_k = p_k (A) r_0 \en where
$p_k$ is a polynomial of degree $k$ satisfying the constraint
$p_k(0) = 1$ that depends on the free coefficients
$\alpha_1, \cdots \alpha_k$.  This is what we would call a
polynomial iteration today.  What comes as a surprise is that
Richardson identifies exactly the problem he has to solve in order
to get a small residual, namely to find a set of $\alpha_i$'s
for which $p_k(t)$ deviates the least from zero, but then does not
invoke Chebysheff's work to find the solution.  He arrives at a
certain solution \emph{`by trial'} for a polynomial of degree~7.
The best solution given knowing that the eigenvalues are in an
interval $[a, b]$ with $a>0$ can be expressed in terms of a Chebyshev
polynomial of the first kind, see, e.g., \cite[Sec. 12.3]{Saad-book2}. 
Chebyshev introduced his polynomials in 1854,
\cite{TchebysheffWorks}, or 56 years prior to Richardson's
article, but his paper addressed completely different issues from
those with  which we are familiar today when analyzing convergence
of certain algorithms or when defining iterative schemes such as the
Chebyshev iteration.
Equally  surprising is the fact that Richardson does not seem to be aware of
the work by
Gauss~\cite{Letter2Gerling} and Seidel~\cite{Seidel1874} on iterative schemes.
His work is truly original in that it defines a completely new method,
the method of polynomial iteration, but misses Chebyshev acceleration
as we know it today.

It was much later that  the missing part was completed
in the work of
Shortley~\cite{Shortley}, Sheldon~\cite{Sheldon},
and finally Golub and Varga~\cite{Golub-Varga}.
This work also lead to a second-order Richardson iteration to accelerate the `basic' iteration
$  u\Sup{k+1} = G u\Sup{k} + f $ which takes the following form: 
\[
  u\Sup{k+1} = \rho \left[ \gamma (G u\Sup{k} + f) + (1-\gamma) u\Sup{k} \right]
  + (1-\rho) u\Sup{k-1}
\]
where, at the difference with the Chebyshev method, the parameters $\rho$ and $\gamma$ are fixed
throughout the iteration.

\section{1930s-1940s: Southwell}
Iterative methods were popularized by a series of papers,
e.g., \cite{GreenSouthwell44,FoxSouthwell45,AllenSouthwell50,SouthwellCyl55,Southwell-1935}, 
and books 
\cite{SouthwellBook,SouthwellBook2} by Richard Southwell and co-workers who put these methods to use
for solving a wide range of problems in mechanical engineering and phsyics.
A good survey of developments with relaxation methods with
a summary of the problems successfully solved by these methods up to the late 1940s
is given by Fox \cite{FoxSurvey48}.
Southwell defined various refinement techniques to standard relaxation, including
block-relaxation (called group relaxation \cite{Southwell-1935}) for example. However, his biggest contribution
was to put these techniques in perspective and to show their effectiveness for handling a large
variety of realistic engineering and physical problems, thus avoiding the use of direct solution methods.
Many of the problems tackled were challenging for that period. 

\section{The SOR era} 
Later toward the mid-20th century the  observation was made that the convergence of a relaxation procedure
could be significantly accelerated by including an \emph{over-relaxation}
parameter. In the language of the iteration \nref{eq:GS}
Over-relaxation (Young and others) is based on the splitting 
\[ \omega A = ( D - \omega E) -( \omega F + (1-\omega )D ) , \]
resulting in a scheme known as the Successive Over-Relaxation (SOR) method.
The 1950s and early 1960s marked a productive 
era for  iterative methods that saw an in-depth 
study of this class of techniques initiated
by David Young and Stanley Frankel.
In a 1950 article, Frankel \cite{Frankel50}   described the
`Liebmann' method, which was just the cyclic relaxation
process described by Seidel, along with an `extrapolated Liebmann method'
which is nothing but the SOR scheme. He discusses the  parameter
$\omega$ (called $\alpha$ in his paper) 
and obtains an optimal value for standard finite difference
discretizations of the Laplacean.
This particular topic received a rather comprehensive treatment 
by David Young in his PhD thesis \cite{Young-thesis} who generalized
Frankel's work to matrices other than those narrowly targetted by Frankel's paper.
The SOR method and its variants, became quite successful, especially
with the advent of digital computing and they enjoyed a
popularity that lasted until the 1980s when Preconditioner Krylov methods
started replacing them.
Here is what Varga \cite{Varga-book}
said about the capabilities of these methods in the year 1960:
\emph{``As an example of the magnitude of problems that have been      
successfully solved by cyclic iterative methods, the Bettis Atomic    
Power Laboratory of the Westinghouse Electric Corporation had in daily     
use in 1960 a 2-dimensional program which would treat as a special case    
Laplacean-type matrix equations of order 20,000."}

He then adds in a  footnote: (paraphrase)   that    
the program was written for the Philco-2000 computer which had    
32,000 words of core storage (32K\-words!) and \emph{``Even more staggering'':}
Bettis had a  3-D code which could treat coupled matrix
equations of order 108,000. This reflects the capability  of iterative methods
and indeed of linear system solvers (direct methods could not handle such systems)
at that point in time.

Up to the early 1980s, this was  \emph{the state of the art} in iterative methods.
These methods  are still  used in some  applications either  as the
main  iterative   solution  method  or  in   combination  with  recent
techniques (e.g.  as smoothers for multigrid or as preconditioners for
Krylov methods).

What I call the SOR era culminated with the production of two major
books that together give a complete view of the state of the art in
iterative methods up to  the late 1960s early 1970s.
The first is by Richard Varga~\cite{Varga-book} which appeared in 1962
and the second by David Young~\cite{Young-book} which appeared in 1971.

\section{A turning point: The Forsythe article}
In  1953,     George   Forsythe published a  great  survey   article
\cite{Forsythe53} in the  Bulleting of the American  Society, with the
title:   \emph{``Solving    linear   algebraic   equations    can   be
  interesting''}.  The paper is rather illuminating as much by the breadth of
its content  as by its  vision.  In  it Forsythe  mentions a  new method,
called the  Conjugate Gradient method,  that appeared on  the horizon.
\emph{  ``  It  is  my  hope,   on  the  other  hand,  to  arouse  the
  mathematician's  interest  by  showing  (sec. 2)  the  diversity  of
  approaches to the solution of (1),  and by mentioning (secs. 3 to 6)
  some problems associated with selected iterative methods.  \emph{The
    newest process on  the roster, the method  of conjugate gradients,
    is outlined  in sec. 7}. Sec.  8 touches on the  difficult general
  subject of errors and "condition," while a few words are hazarded in
  sec. 9  about the effect  of machines  on methods.''}  The  title is
intriguing but what it  is even more so when one  reads the comment by
the  author   that  the   title  of   the  submitted   manuscript  was
\emph{``Solving linear  systems is  not trivial''}.  We  will probably
never know the reason for the change, but it seems clear that in those
days,  solving linear  systems of  equations  could be  thought to  be
`trivial' from some angle.  \footnote{When I was working for my PhD in
  France, I  was once asked  about the topic of  my thesis and  when I
  replied that it  was about sparse matrix methods I  was told ``..but
  the problem of solving linear  systems of equations is solved. Isn't
  that  just  tinkering?''   Just   like  Young  in  \cite{YoungHist},
  ``though discouraged I continued to work.''}

George Forsthye (1917-1972) joined Stanford in 1957 (Math) and founded the computer science department,
one of the first in the nation, in 1965. Knuth discusses this era in
\cite{Knuth72} illustrating the fact that it was not an easy task to start a computer science department at the time
and praising Forsythe's vision. George Forsythe is now considered the father of modern numerical analysis.

\section{In Brief: Chaotic Iterations} 
In the early days of electronic computing many people started to see the potential of parallel processing.
\emph{Chaotic relaxation} was viewed as a way of exploiting this avenue.
It is interesting to see how early this vision of parallelism emerged.
Two papers introduced the term \emph{``free steering''} for a relaxation method in which the components to
be relaxed are chosen freely one by Ostrowski in 1955 \cite{Ostrowski55}  and the other by Schechter
in 1959 \cite{Schechter59}. Both studied convergence for H-matrices.
The article by Chazan and Miranker in 1969 introduced the term \emph{``chaotic relaxation''} which was
adopted for a while until it was replaced later by the term \emph{``asynchronous relaxation''}.

Here is a quote from  the paper by Chazan and Miranker [1969] that explains the motivation and context of their work:
\emph{``The problem of chaotic relaxation was suggested to the authors by
J. Rosenfeld who has conducted extensive numerical experiments with
chaotic relaxation} [J. Rosenfeld (1967)].
\emph{Rosenfeld found that the use of more processors
decreased the computation time for solving the linear system.
(...)
The chaotic form of the relaxation eliminated
considerable programming and computer time in coordinating the processors
and the algorithm. His experiments also exhibited a diminution
in the amount of overrelaxation allowable for convergence in a chaotic mode."}

The  article \cite{Rosenfeld67} by Rosenfeld, mentioned above, 
 seems  to be  the  first to  actually
implement  chaotic  iteration  on  a parallel  machine  and  show  its
effectiveness and potential.   A group from the  French school started
this line of  work with the doctoral thesis of  Jean-Claude Miellou in
1970. 
Miellou had a series of articles  in the Comptes Rendues de l'Academie
des Sciences  (Proceedings of  the French  Academy of  Sciences), see,
e.g.,  \cite{MiellouCras1,MiellouCras2}.   The paper  \cite{Miellou75}
summarizes some of the work he  did on chaotic relaxation methods with
delay.   The  work by  Miellou  generated  an important  following  in
France, with  papers that focused  on convergence as well  as parallel
implementations     \cite{Baudet78,RobertCharnayMusy,Eltarazi,ElBaz90,
  ElbazBertsekas,SpiteriGiraud,BahiMiellou93,MiellouEtAl98}.   Some of
the  work  done  in  France   in  those  days  was  truly  visionnary.
Discussions that I  attended as a student in Grenoble,  could be tense
sometimes, with one camp claiming  that the methods were utopian. They
were not  necessarily utopian but certainly  ahead of their time  by a
few decades.   In fact this  work has  staged a  strong come
back with  the advent of  very large high-performance  computers where
communication is expensive, see, for example,
\cite{Magoules,FrommerSzyld,EmadAl,AntzAl,FrommerAsyn,ChowAlAsyn} among many others.

\section{Meanwhile, on the opposite camp}
At this point, the reader may be led to believe that direct methods,
or direct elimination methods using Gauss' terminology, were about to be  abandoned 
as the success of iterative methods was spreading to more areas
of engineering and science. In fact, the opposite happened.
Developers of direct methods for sparse linear systems became very active
starting in the 1960s and the whole area witnessed an amazing progression 
in the few decades that followed.
Here, it is  good to begin  by mentioning the survey  article by Iain Duff
\cite{Duff-survey} which had over 600 references already in  1977.
Major advances were made for general sparse matrices -- as opposed to those
matrices with regular structure that came from finite difference techniques applied to PDEs on
simple regions.
When a sparse linear system is solved by Gaussian elimination, some zero entries
will become nonzero and because of the repetitive linear combination of rows
the final matrix may loose sparsity completely. A new non-zero entry created by the process
is called a `fill-in' and the number of fill-ins created depends enormously on the way the equations
and rows are ordered. Then, a big part of the know-how in \emph{`sparse direct methods'} 
is to try to find orderings that minimize fill-in.

The discovery of sparse matrix techniques began with the link made between graph
theory and sparse Gaussian elimination by Seymour Parter \cite{Parter} in 1961.
This paper gave a model of the creation of fill-in that provided a better understanding
of the process. 
Graphs played a major role thereafter but 
it took some time before a major push was made to
exploit this link in the form of a theorem that guarantees the non-creation
of fill-ins by judicious reordering~\cite{RoseTarjanLueker,Rose-Tarjan}.
One important feature of sparse direct methods that distinguishes them from  iterative methods,
is that they are rather complex to  implement. Today, it takes man-years of effort to develop a good working code with
all the optimized features that have been gathered over years of steady progress. In contrast, it would
take a specialist a few days or weeks of work to develop a small set of preconditioners (e.g., of ILU-type) with
one or two accelerators.  This distinction has had an impact on available software.
In particular,  sparse direct solvers (SPARSPAK, YSMP, ..) were all commercial packages at the beginning.

A major contribution, and boost to the field, was made in 1981 by Alan George and Joseph Liu who published
an outstanding book  \cite{George-Liu-Book} that layed out all that has been learned
on sparse direct solution techniques for solving symmetric linear systems up to that point. The book also  included FORTRAN routines
 and this lead to the first package,
 called SPARSPAK \cite{SPARSPAK}, for solving sparse symmetric positive definite linear systems\footnote{
   As was just mentioned SPARSPACK  was a commercial package but the book included listings of the main routines.}.

The speed with which progress was made at the early stages of research on sparse direct solvers is staggering.
Table~\ref{tab:Liu}, reproduced from \cite{George-Liu-evol}, shows the evolution of the performance of the minimum degree algorithm,
a reordering technique to reduce fill-in in Gaussian elimination from its inception to 1989.
With each discovery, or new trick, there is a gain, often quite substantial, in performance, both in the reduction of the number of
nonzero entries and the time of the procedure.  Since 1989, many more new ingredients have been found that make sparse direct solvers
hard to beat  for certain types of problems.

\begin{table}[h!]
\begin{center} 
  \begin{tabular}{|l|l|r|r|} \hline 
  Version  &    Minimum Degree Algorithm  &       Off-diagonal& Ordering\\
           &                               &   Factor Nonz    &Time\\  \hline 
Md\#1 & Final minimum degree                 &   1,180,771  &   43.90\\
Md\#2 &   Md\# 1 without multiple elimination    &   1,374,837  &   57.38\\
Md\#3 &   Md\# 2 without element absorption      &   1,374,837  &   56.00\\
Md\#4 &   Md\# 3 without incomplete deg update   &   1,374,837  &   83.26\\
Md\#5 &   Md\# 4 without indistinguishable nodes &   1,307,969  &  183.26\\
Md\#6 &   Md\# 5 without mass elimination        &   1,307,969  & 2289.44\\ \hline 
\end{tabular}
\end{center} 
\caption{Evolution of the minimum degree algorithm up to 1989 according to~\cite{George-Liu-evol}. 
  \label{tab:Liu}}
\end{table} 

The merits and disadvantages of direct and iterative methods have  been compared  since the earliest
paper of Gauss, see Section~\ref{sec:Gauss}.
In his 1959 book \cite[p. 173]{Bodewig-Book} Bodewig states that
\emph{``Compared with direct methods, iteration methods have the great
  disadvantage that, nearly always, they converge too slowly and, therefore, the number of operations is large''}. Then he continues that
in fact
\emph{``For most systems the iteration does not converge at all. The methods for making convergent an
  arbitrary system are circumstantial.''}

The only potential advantage of iterative methods  over direct methods he saw was that
\emph{``Rounding errors cannot accumulate, for they are restricted to the
last operation. So, without doubt, this is an advantage compared with
direct methods. Yet this advantage costs probably more than it is worth.''}

Later, David M.  Young \cite{Young-book} states in the  first chapter of
his book (1971): \emph{The use of direct methods even for solving very
  large problems  has received  increased attention recently  (see for
  example Angel, 1970). In some  cases their use is quite appropriate.
  However, \emph{there  is the  danger that if  one does  not properly
    apply  iterative  methods  in  some  cases  one  will  incorrectly
    conclude that they are not effective} and that direct methods must
  be used.   It is  hoped that  this book  will provide  some guidance
  (...) } A  comparison from the opposite camp (George  \& Liu's book
\cite{George-Liu-Book})
warns that:  \emph{ (...) 
  Unless the  question of  which class  of methods  should be  used is
  posed in a quite narrow and  well defined context, it is either very
  complicated  or  impossible  to   answer.}  The  authors  then  give
reference to Varga and Young and say that there are no books on direct
solvers and \emph{ ``In addition, there are situations where it can be
  shown quite convincingly that direct  methods are far more desirable
  than any conceivable iterative scheme.''}
Surprisingly, this section of the book does not mention
the relative ineffectiveness  of direct solvers for  large 3D problems
(though this was clearly known by the authors at the time, see below).

The debate  has somewhat diminished  recently with the  consensus that
iterative methods are competitive for 3-D problems -- but that for 2-D
problems the benefits  may be outweighed by the lack  of robustness of
these methods  for indefinite problems.  The  common argument that
is given to  illustrate this fact is  to compare the result  of one of
the best  orderings for regular  grids in the  2-D and 3-D  cases,
as illustrated in  \cite{George-Liu-Book}.   Consider a  standard
Poisson equation on an  $n \times n$ regular grid in  2-D, and then on
an $n \times n \times n$ regular grid  in 3-D. We call $N$ the size of
the resulting system,  so $N = n^d$ where $d$  is the space dimension,
i.e., $d=  2, 3$.   The order of  the cost is  given by  the following
table:
\begin{center} 
\begin{tabular}{l|c|c} 
                  &  2-D            &  3-D \\ \hline 
space (fill)  \ \   & \ { $ O(N \log N)$ } \ 
 & \ {  $O(N^{4/3}) $ }  \  \\  \hline  
time  (flops) \ \   & {  \ $ O(N^{3/2})  \ $ }
 & {   \ $O(N^{2})   $ } \  \\ 
\end{tabular} 
\end{center} 
The table shows a  significant  difference in complexity between the 2-D and the 3-D cases.

A widespread  misconception is that \emph{3-D problems are harder just
  because they are bigger.} In fact they are just \emph{intrinsically} harder as is suggested in the above table.
When I teach sparse matrix techniques in a numerical linear algebra course,
I often give a demonstration in Matlab to illustrate direct solution methods. I show a live illustration of using
the \emph{back-slash} operation \footnote{In Matlab  a sparse linear system $Ax=b$ can be solved by the command
  $x = A \backslash b $. This back-slash operation will invoke a sparse direct solver to produce the answer.} 
in matlab to solve a linear system involving a coefficient matrix that comes from a centered difference
discretization of Poisson's equation on a 2-D or a 3-D mesh. The idea is to show that for the same size problem,
e.g., $350 \times 350$ grid in 2-D versus $50 \times 50 \times 49$ grid in 3-D (leading to a problem of size $N=122,500$ in each case),
the 3-D problem takes much longer to solve. For this example it can take ~11 sec. for the 3-D problem and ~ 0.7 sec for the
2-D problems on my laptop. What I also tell the audience is that in past years I was gradually increasing the size of these
problems as times went down. A decade ago for example, I would  have a demo with a problem of size approximately $N=20,000$ if I wanted
not to have students wait too long for the answer. Of course, this gain reflects progress in both hardware and algorithms.

\section{One-dimensional projection processes}\label{sec:Proj1D}
The method of  steepest descent was introduced by Cauchy  in 1829 as a
means of  solving a  nonlinear equation  related to  a problem  of the
approximation of an integral.  A detailed account of the origin of the
steepest descent  method is  given in \cite{PetrovaOnCauchy}  where we
learn that  Riemann, Nekrasov,  and later  Debye were  also associated
with the  method.  In 1945  Kantorovitch introduced the method  in the
form we know today for linear systems for SPD matrices: \eq{eq:Kantor}
\min_x \ f(x) \equiv  \frac{1}{2} x^T A x - b^T x  .  \en The gradient
of  the above  function is  $ \nabla  f(x) =  A x  - b$  which is  the
negative of the residual $b - Ax$ and so the steepest method will just
generate an iteration of the form
\[x_{k+1} = x_k + \omega_k r_k \]
where $r_k = b - Ax_k$ and $\omega_k$ is selected to minimize \nref{eq:Kantor} at each step. Convergence can easily shown for matrices that are
symmetric positive definite. Methods of this type are one-dimensional
projection methods in the sense that they produce a new iterate
$x\Sup{new}$
from a current iterate $x$ by a modification of the form
$x\Sup{new} = x + \delta$ where $\delta $ belongs to a subspace of
dimension $1$. In the case of the steepest descent method
we can write $\delta = \alpha \nabla f(x) = \alpha (b - A x)$
and it is easy to calculate $\alpha $ if we wish to minimize
\nref{eq:Kantor}. 

Simple projection methods of this type for solving linear systems were proposed earlier.
For example, in a short paper \cite{Kaczmarz}, Kaczmarz
described a method in which
at each step $d$ is selected to be the vector $a_i = A^T e_i$, the
$i$-th row of $A$ written as a column vector. In this case,
\eq{eq:Kacz}
x\Sup{new} = x + \alpha a_i \quad \alpha = \frac{r_i}{\| a_i \|_2^2} 
 \en 
 in which $r_i $ is the $i$-th component of the current residual vector
 $b - A x$. Equation is written in a form that avoids clutter but
 we note that the indices $i$ of the components that are modified are
 cycled from 1 to $n$ and this is repeated until convergence.
 We can rewrite \nref{eq:Kacz} as
 $ x\Sup{new} = x + \frac{ e_i^T r }{ \| A^T e_i \|_2^2 }
 A^T e_i$. If $x_*$ is the exact solution and we write
 $x\Sup{new} = x + \alpha a_i$ then we have
 \begin{align}
 \| x_* - x\Sup{new} \|_2 ^2
 & = \| (x_* - x) - \alpha a_i \|_2^2   \nonumber \\
   & = \| x_* - x\|_2^2
     + \alpha^2 \| a_i \|_2^2 -  2 \alpha < x_* - x ,  a_i > . 
\label{eq:kaczE} 
 \end{align}
 This is a quadratic function of $\alpha$ and the minimum is reached when
 \[
   \alpha = \frac{ < x_* - x ,  a_i >}{\| a_i \|^2} 
   = \frac{< x_* - x ,  A^T e_i >}{\| a_i \|^2} 
   = \frac{ < A (x_* - x) , e_i >}{\| a_i \|^2}
   =   \frac{ < r, e_i >}{\| a_i \|^2} \]
 which is the choice made in the algorithm. In addition, Kaczmarz was
 able to show convergence. Indeed, with the optimal $\alpha$,
 equation \nref{eq:kaczE} yields
 \eq{eq:kacz2} \| x_* - x\Sup{new} \|_2 ^2 =  \| x_* - x\|_2^2 -
 \frac{r_i^2 }{\| a_i \|_2^2 } ,
\en 
showing that the error must decrease. From here there are a number of ways
of showing convergence. The simplest is to observe that the norm of the error
$\|x_* - x\Sup{j} \|$ must have a limit and therefore \nref{eq:kacz2} implies
that each residual component $r_i$ converges to
zero, which in turn implies that $x\Sup{j}$ converges to the
solution.
The method is motivated by a simple interpretation.
The solution $x$ is located at the intersection of the
$n$ hyperplanes represented by the equations $ a_i x - b_i = 0$ and 
the algorithm projects the current iterate on one of these hyperplanes
in succession, bringing it closer to the solution each time.

At almost the same time, in 1938,  Cimmino 
\cite{Cimmino} proposed a  one-dimensional process which has some
similarity with the Kaczmarz algorithm. Instead of projecting the
solution onto the various hyperplaces, Cimmino generates $n$ intermediate
solutions each of which is a mirror image of the current iterate with
respect to the hyperplanes. Once these are available then he takes
their convex combination. 
Specifically,
Cimmino defines intermediate iterates in the form 
\begin{equation} 
x\Sup{j} = x + 2 r_j a_{j,.}  
\end{equation}
where $r_j$ is the $j$-th component of the residual $r = b - Ax$,
and then takes as a new iterate a convex combination of these points:
\begin{equation} 
x_{new}  = \sum \mu_j x\Sup{j} . 
\end{equation}
Details on this method and on the life and contributions of
Cimmino can be found in Michele Benzi's article~\cite{BenziOnCimmino}.

\section{Krylov methods take off: The CG algorithm}
One-dimensional projection methods and Richardson
iteration are of the form $ x_{k+1} = x_k + \beta_k d_k$, where $d_k$ is a
certain direction that is generated from the current iterate only.

It was Frankel  who in 1950 had the idea to extend these
to  a `second-order' iteration of the form \cite{Frankel50} 
\eq{eq:frankel} x_{k+1} = x_k + \beta_k d_k \quad \mbox{where}\quad 
d_k = r_k - \alpha_k d_{k-1}.
\en
Frankel was inspired by the solution of time-dependent partial differential equations
such as the heat equation which allowed him to add a parameter.
We can recover the Chebyshev iteration by using constant coefficients
$\alpha_k $ and $\beta_k$ as we saw before.
A method of the type represented by \nref{eq:frankel} with constant coefficients
$\alpha_k$, $\beta_k$ was termed \emph{semi-iterative} method.
In his 1957 article Varga \cite{Varga57}
uses this term for any polynomial method
and mentions earlier work by Lanczos and Stiefel.
The 1961 paper by Golub and Varga \cite{Golub-Varga} explains how Chebyshev
polynomials can be used effectively and stably.

The understanding and development of semi-iterative methods is deeply rooted in approximation theory.
The residual of the approximation $x_{k+1}$ obtained from  Richardson type iteration of the form
$x_{k+1} = x_k + \omega_k r_k$, can be shown to be equal to
\[
  r_{k+1} = (I -\omega_k A) (I -\omega_{k-1} A) \cdots  (I -\omega_{0} A) r_0
  \equiv p_{k+1} (A) r_0
\]
where  $ p_{k+1}$  is  a  polynomial of  degree  $k+1$ satisfying  the
conditin  $p_{k+1} (0)  =  1  $. One  can  therefore design  effective
iterative schemes by selecting polynomials of this type that are small
on a set that contains the  spectrum of $A$.  Many papers adopted this
`approximation theory'  viewpoint. This  is most apparent  in Lanczos'
work. Thus,  the remarkable 3-term  recurrence obtained by  Lanczos to
generate an orthogonal  basis of the Krylov subspace  is a consequence
of the Stieljes procedure for generaing orthogonal polynomials.
Magnus Hestenes [UCLA] and Eduard Stiefel [ETH,  Z\"urich]
developed the method of Conjugate Gradient independently.
The article~\cite{Saad-vdv-reviewY2K} describes how the two authors
discovered that they both developed the exact same method independently
at the occasion of a conference held at UCLA in 1951.
Lanczos developed another method that exploited what we now call Lanczos
vectors, to obtain the solution from the Krylov subspace that has the smallest
residual norm. His paper \cite{Lanczos52} appeared within 6 months of the one by
Hesteness and Stiefel. 
The method developed by Lanczos is mathematically equivalent to
what we would call the Minimal Residual method today, but it is
implemented with the Lanczos procedure.

Though  not perceived  this way  at the  time, the  conjugate gradient
method was the  single most important advance made in  the 1950s.  One
of the main  issues with  Chebyshev semi-iterative  methods is  that they
require fairly accurate estimates of extremal eigenvalues, since these
define the interval in which the residual polynomial is minimized.
The conjugate gradient method bypassed this drawback -- but it was viewed
as  an unstable, direct method.
Engeli et  al. \cite{engirust59} were the first to view  the method as
an iterative process and indicated that this process 
can take $2n$ to $3n$ steps to 'converge'.

The method laid dormant until the early 1970s when a  paper
by John Reid   \cite{Reid71}
showed the practical interest of this iterative viewpoint
when considering large sparse linear systems.
With the advent of  incomplete Cholesky preconditioners
developed by Meijerink and van der Vorst in 1977
\cite{Mei-vdv}, the method gained tremendous popularity and
ICCG (Incomplete Cholesky Conjugate Gradient) became
de facto iterative solver for the general Symmetric Positive Definite case.

\section{Krylov methods: the `nonsymmetric' period}
Nonsymmetric linear systems were given less attention right from the beginning
of iteative methods.
In his 1952 paper,   Lanczos \cite{Lanczos52} discusses  a method that
is essentially equivalent to what we now call the BiCG algorithm  and
then drops the method by stating:
\emph{let us restrict our attention to symmetric case ... (Normal equations.)}.

However, the demand for nonsymmetric solvers started to strengthen when 
applications in aerospace engineering for example gained
in importance. Thus, the success of  the CG method lead researchers to
investigate Krylov subspace methods  for the nonsymmetric case.

It was only in  1976 that Fletcher introduced the \emph{BiCG}, which was
based on the Lanczos process. 
The BiCG uses two matrix-vector products: one with $A$ and the other with
$A^T$. However,
the operations with $A^T$ are only needed to generate the coefficients
needed for the projection and they were therefore viewed as wasteful.
A number of methods later appeared whose goal was to avoid these products.
The first of these  was the Conjugate Gradient Squared (CGS) \cite{sonn} 
developed by Sonneveld in 1984. Then came the BiCGSTAB 
\cite{vorstbicg} in 1992, along with variants, e.g.,
\cite{slfo93},  QMR \cite{Freund-Nachtigal-QMR},
TFQMR \cite{Freund-TFQMR}, QMRSTAB,  \cite{chgasi94} and several others.
These contributions to accelerators for the nonsymmetric case
are  described in detail in the earlier paper ~\cite{Saad-vdv-reviewY2K} which covers the period
of the twentieth century. In fact research on accelerators has  been less active since 2000
while preconditioners have attracted continued attention.

\section{Present and future}
Modern numerical linear  algebra started with the  influence of George
Forsythe and one could view  his 1953 survey paper \cite{Forsythe53}
 as  a sort  of road-map.  Since  then, the  field has
changed directions several times, often  to respond to new demands from
applications.  So the  natural question  to  ask is  \emph{`what next?'}   For
iterative  methods,   research  is  still   active  in  the   area  of
preconditioners  for  some  types  of  problems  (Helmholtz,  Maxwell,
Structures,...), as well in  developing efficient parallel algorithms.
For example,  it was  noted earlier  that asynchronous  iterations are
back.   On  the other  hand  research  on accelerators  has  subsided.
Another  observation is  that  the fields  of  numerical analysis  and
numerical  Linear Algebra  are gradually  disapperaring from  computer
science graduate programs.  This is unfortunate because some research topics
fit  better in  computer  science than  in  mathematics.  Among  these
topics  we can  mention: sparse  matrix techniques  and sparse  direct
solvers,   preconditioning   methods,  effective   parallel   solvers,
and graph-based  methods.  Some  of  these topics  may  reappear in  other
areas,  e.g.  computational  statistics,  and  machine learning, but it
they are  no  longer  represented in  either  computer science  or
mathematics, there will be a lack of  students trained in them.

When trying to answer the questions  ``What next?'' we need to remember
that for the bigger part of  the 20th century, solution
techniques (iterative and direct) were aimed primarily at solving certain types of
PDEs, and this
was driven  in part by  demand in some engineering  applications, most
notably  the aerospace,  the  automobile, and the semi-conductor
industries.   Therefore, a related question to ask  is ``What
new demands  are showing up  at the  horizon?''.  
Currently, the answer to this question is without a doubt related to the
emergence of data mining and machine learning. 
Conferences that used to bear the title
`computational X'   in the past are now often replacing this title,
or augmenting it, by
`machine learning X'. Linear algebra is gradually addressing
tasks that arise in the optimization, and computational statistics
problems of machine learning. The new linear algebra specialist 
encounters such problems as evaluating matrix functions,
computing and updating the
SVD, fast low-rank approximation methods, random sampling methods, etc.
An important  new consideration in all of these topics is the
pre-eminence of randomness and stochastic approaches.
In this context, methods such as the conjugate gradient or GMRES,
that are based on global optimality are not  adapted to randomness
and it may be time to look for alternatives or to reformulate them.
There are opportunities also in adapting various techniques learned
in linear algebra and more broadly in numerical analysis to solve
various problems in machine learning. One can
echo the title of Forsythe's 1953 paper  \cite{Forsythe53} by saying that
\emph{``Solving matrix problems in machine learning can be interesting''.} 

    \nocite{Saad-vdv-reviewY2K,Golub-Oleary,Axelsson-survey,Forsythe53}

\nocite{BenziOnCimmino}

\bibliographystyle{siam} 
\bibliography{/home/saad/tex/BIB/strings1,/home/saad/tex/BIB/bib,/home/saad/tex/BIB/newbib,/home/saad/tex/BIB/saad,/home/saad/tex/BIB/henka,/home/saad/tex/BIB/henkb,/home/saad/tex/BIB/henkc,/home/saad/tex/BIB/sparse,/home/saad/tex/BIB/ilu,/home/saad/tex/BIB/review2K,local}

\end{document}